\documentclass{amsart}
\usepackage[section]{placeins}
\usepackage{float}
\usepackage{tcolorbox}
\usepackage{enumerate}
\usepackage{mathtools,amssymb}
\usepackage{subfig}
\usepackage{dcolumn}
\usepackage{multirow}
\usepackage{tikz-cd}
\usepackage[all]{xy}
\usepackage{graphicx}
\usepackage{hyperref}

\usepackage{easyReview}

\newtheorem{theorem}{Theorem}[section]

\newcommand{\Vect}{{\mathsf{Vect}}}
\newcommand{\Proj}{{\mathsf{Proj}}}
\newcommand{\Mod}{{\mathsf{Mod}}}
\newcommand{\Coh}{{\mathsf{Coh}}}
\newcommand{\Perf}{{\mathsf{Perf}}}
\newcommand{\MF}{{\mathsf{MF}}}
\newcommand{\Fuk}{{\mathsf{Fuk}}}
\newcommand{\Singu}{{\mathsf{Sing}}}

\newcommand{\XX}{{\mathcal{X}}}
\newcommand{\EE}{{\mathcal{E}}}
\newcommand{\FF}{{\mathcal{F}}}
\newcommand{\OO}{{\mathcal{O}}}
\newcommand{\sslash}{\mathbin{/\mkern-6mu/}}

\newcommand{\DR}{{\mathsf{DR}}}
\newcommand{\Res}{{\mathsf{Res}}}
\newcommand{\Rep}{{\mathsf{Rep}}}
\newcommand{\Tw}{{\mathsf{Tw}}}
\newcommand{\Jac}{{\mathsf{Jac}}}
\newcommand{\Hom}{{\mathsf{Hom}}}
\newcommand{\End}{{\mathsf{End}}}
\newcommand{\Sp}{{\mathsf{Sp}}}
\newcommand{\orb}{{\mathsf{orb}}}
\newcommand{\age}{{\mathsf{age}}}
\newcommand{\SL}{{\mathsf{SL}}}
\newcommand{\SU}{{\mathsf{SU}}}
\newcommand{\SO}{{\mathsf{SO}}}
\newcommand{\Spec}{{\mathsf{Spec}}}

\newcommand{\aff}{{\mathbb{A}}}

\newcommand{\complex}{{\mathbb{C}}}
\newcommand{\integers}{{\mathbb{Z}}}
\newcommand{\rationals}{{\mathbb{Q}}}
\newcommand{\naturals}{{\mathbb{N}}}

\title{Revisiting the Classical McKay Correspondence, Derived Equivalences and the Spectrum of Kleinian Surface Singularities: A Look Through the Mirror}

\author{Enrique Becerra}

\author{Ludmil Katzarkov}

\author{Ernesto Lupercio}

\begin{document}

\keywords{Homological mirror symmetry, Mckay correspondence, spectrum.}
\subjclass[2020]{14-02,14E16,14F08, 14J33,14B05}
	
	\maketitle
	
	\begin{center}{\footnotesize{\emph{To Kaxjuu, on the occasion of her first nine months in this world.}}}\end{center}
	
	\begin{abstract} In this article, we revisit the classical McKay correspondence via  homological mirror symmetry. Specifically, we demonstrate how this correspondence can be articulated as a derived equivalence between the category of vanishing cycles associated with a Kleinian surface singularity and the category of perfect complexes on the corresponding quotient orbifold. We further illustrate how this equivalence allows for the interpretation of the spectrum of a Kleinian surface singularity solely in terms of the representation-theoretic data of the associated binary polyhedral group.
	\end{abstract}

	\section{Introductory Overview}
	
	\subsection{Kleinian Surface Singularities}
	
	In his seminal 1884 work \emph{``Lectures on the Icosahedron and Equations of the Fifth Degree''} \cite{klein1884vorlesungen}, F. Klein accomplished, among many other things, the classification of all finite subgroups of $\SL(2, \mathbb{C})$. Klein demonstrated that every finite subgroup $G$ of $\SL(2, \mathbb{C})$ must be isomorphic to one of three types: a cyclic group, a binary dihedral group, or a binary polyhedral group.

	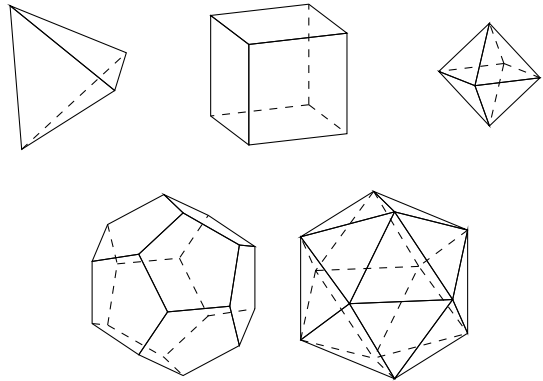
\begin{figure}[H]
		\begin{tikzpicture}[scale=0.7]
			
			\coordinate (T0) at (-5.10427589776, 1.3096936302264104);
			\coordinate (T1) at (-3.11649858992, -0.3031514993551603);
			\coordinate (T2) at (-4.88350141008
			, -1.4214862452202075);
			\coordinate (T3) at (-2.89572410224
			, 0.41494411434895745);
			
			\draw (T0) -- (T1) -- (T2) -- cycle;
			
			\draw (T0) -- (T3) -- (T1) -- cycle;
			
			\draw[dashed] (T2) -- (T3);
			
			\coordinate (H0) at (-0.5696813314905527, 0.5728160673358398);
			\coordinate (H1) at (-0.5696813314905527, -1.3378569109153722);
			\coordinate (H2) at (1.2943968404439, 0.7869841443124108);
			\coordinate (H3) at (1.2943968404439, -1.1236888339388011);
			\coordinate (H4) at (-1.2943968404439, 1.1236888339388011);
			\coordinate (H5) at (-1.2943968404439, -0.7869841443124108);
			\coordinate (H6) at (0.5696813314905527, 1.3378569109153722);
			\coordinate (H7) at (0.5696813314905527, -0.5728160673358398);
			
			\draw (H0) -- (H1) -- (H3) -- (H2) -- cycle;
			
			\draw (H0) -- (H1) -- (H5) -- (H4) -- cycle;
			
			\draw (H0) -- (H2) -- (H6) -- (H4) -- cycle;
			
			\draw[dashed] (H7) -- (H6);
			\draw[dashed] (H7) -- (H5);
			\draw[dashed] (H7) -- (H3);
			
			\coordinate (O0) at (3.72483666195, -0.21095489306059886);
			\coordinate (O1) at (4.27516333805159693, 0.21095489306059886);
			\coordinate (O2) at (3.03860250812, 0.06037778654840108);
			\coordinate (O3) at (4.9613974918795568, -0.06037778654840108);
			\coordinate (O4) at (4.0, -0.9756293127952373);
			\coordinate (O5) at (4.0, 0.9756293127952373);
			
			\draw (O4) -- (O2) -- (O0) -- cycle;
			
			\draw (O4) -- (O3) -- (O0) -- cycle;
			
			\draw (O5) -- (O2) -- (O0) -- cycle;
			
			\draw(O5) -- (O3) -- (O0) -- cycle;
			
			\draw[dashed] (O1) -- (O2);
			\draw[dashed] (O1) -- (O3);
			\draw[dashed] (O1) -- (O4);
			\draw[dashed] (O1) -- (O5);
			
			\coordinate (D0) at (-0.74914330421, -3.23965255503);
			\coordinate (D1) at (-0.74914330421, -5.1503255332779627);
			\coordinate (D2) at (-1.34018371754, -2.67501008163);
			\coordinate (D3) at (-1.34018371754, -4.5856830598829275);
			\coordinate (D4) at (-2.6598162824642664, -3.41431694012);
			\coordinate (D5) at (-2.6598162824642664, -5.3249899183682845);
			\coordinate (D6) at (-3.2508566957869456, -2.84967446672);
			\coordinate (D7) at (-3.2508566957869456, -4.7603474449732492);
			\coordinate (D8) at (-1.81735846792, -2.62871720992);
			\coordinate (D9) at (-1.81735846792, -5.7202510301231948);
			\coordinate (D10) at (-2.18264153207910092, -2.27974896987);
			\coordinate (D11) at (-2.18264153207910092, -5.3712827900732547);
			\coordinate (D12) at (-0.93140784028, -4.40283109341752804);
			\coordinate (D13) at (-1.88773131777, -3.48922037997);
			\coordinate (D14) at (-2.11226868223217829, -4.5107796200274473);
			\coordinate (D15) at (-3.0685921597130592, -3.59716890658);
			\coordinate (D16) at (-0.4542330899, -3.26826312318);
			\coordinate (D17) at (-3.5457669100982248, -3.55087603488);
			\coordinate (D18) at (-0.4542330899, -4.44912396512249836);
			\coordinate (D19) at (-3.5457669100982248, -4.7317368768227392);
			
			\draw (D0) -- (D8) -- (D10) -- (D2) -- (D16) -- cycle;
			
			\draw (D4) -- (D8) -- (D10) -- (D6) -- (D17) -- cycle;
			
			\draw (D0) -- (D8) -- (D4) -- (D14) -- (D12) -- cycle;
			
			\draw (D16) -- (D0) -- (D12) -- (D1) -- (D18) -- cycle;
			
			\draw (D1) -- (D12) -- (D14) -- (D5) -- (D9) -- cycle;
			
			\draw (D4) -- (D14) -- (D5) -- (D19) -- (D17) -- cycle;
			
			\draw[dashed] (D19) -- (D7) -- (D11) -- (D9);
			\draw[dashed] (D7) -- (D15) -- (D6);
			\draw[dashed] (D15) -- (D13) -- (D2);
			\draw[dashed] (D13) -- (D3) -- (D18);
			\draw[dashed] (D3) -- (D11);
			
			\coordinate (I0) at (2.19866933079506122, -2.60892813697);
			\coordinate (I1) at (2.19866933079506122, -5.78049020534);
			\coordinate (I2) at (1.8013306692, -2.21950979466);
			\coordinate (I3) at (1.8013306692, -5.39107186303);
			\coordinate (I4) at (3.301520307589847, -4.27557655385);
			\coordinate (I5) at (2.6586128480926363, -3.64548444015);
			\coordinate (I6) at (1.34138715191, -4.35451555985);
			\coordinate (I7) at (0.69847969241, -3.72442344615);
			\coordinate (I8) at (3.5857810341849234, -2.95607042479);
			\coordinate (I9) at (0.41421896581, -3.08379641953);
			\coordinate (I10) at (3.5857810341849234, -4.91620358047);
			\coordinate (I11) at (0.41421896581, -5.04392957521);
			
			\draw (I0) -- (I2) -- (I8) -- cycle;
			\draw (I0) -- (I2) -- (I9) -- cycle;
			\draw (I0) -- (I8) -- (I4) -- cycle;
			\draw (I0) -- (I9) -- (I6) -- cycle;
			\draw (I0) -- (I4) -- (I6) -- cycle;
			
			\draw (I9) -- (I11) -- (I6) -- cycle;
			\draw (I1) -- (I11) -- (I6) -- cycle;
			\draw (I1) -- (I4) -- (I6) -- cycle;
			
			\draw (I4) -- (I10) -- (I8) -- cycle;
			\draw (I4) -- (I10) -- (I1) -- cycle;
			
			\draw[dashed] (I1) -- (I3);
			\draw[dashed] (I11) -- (I3) -- (I10);
			\draw[dashed] (I11) -- (I7) -- (I9);
			\draw[dashed] (I2) -- (I7) -- (I3);
			\draw[dashed] (I5) -- (I7);
			\draw[dashed] (I2) -- (I5) -- (I3);
			\draw[dashed] (I8) -- (I5) -- (I10);
			
		\end{tikzpicture}
		\caption{The five Platonic solids.}
		\label{fig:Platonic}
	\end{figure}
	
	We can arrive at this classification as follows. Choose a $G$-invariant Hermitian metric on $\mathbb{C}^2$. Consequently, any such $G$ becomes conjugate to a subgroup of the corresponding special unitary group $\SU(2) \subset \SL(2, \mathbb{C})$. On the other hand, we have the group extension:
	\begin{center}
		\begin{tikzcd}
			1 \ar[r] & \mathbb{Z}/2\mathbb{Z} \ar[r] & \SU(2) \ar[r, "\pi"] & \SO(3) \ar[r] & 1
		\end{tikzcd}
	\end{center}
	
	This extension can be elegantly realized through Hamilton's algebra of quaternions:
	$$\mathbb{H} \cong \left\{ q = \begin{pmatrix}
		z & w\\
		-\overline{w} & \overline{z}
	\end{pmatrix} : z, w \in \mathbb{C} \right\} \subset M_2(\mathbb{C})$$
	
	Note that one can identify the special unitary group $\SU(2)$ with the group of unit quaternions $\mathbb{H}^* := \left\{ q \in \mathbb{H} : \|q\|^2 := \det(q) = 1 \right\} \cong \SU(2)$, which acts on the Euclidean space $\mathbb{R}^3$ by conjugation after identifying the latter with the space of purely imaginary quaternions:
	$$\mathbb{R}^3 \cong \left\{ \begin{pmatrix}
		r\sqrt{-1} & z\\
		-\overline{z} & -r\sqrt{-1}
	\end{pmatrix} : r \in \mathbb{R}, z \in \mathbb{C} \right\} \subset \mathbb{H}.$$
	
	Thus, the assignment of taking a unit quaternion $q \in \mathbb{H}^*$ to the corresponding conjugation operator $v \in \mathbb{R}^3 \mapsto q \cdot v \cdot q^{-1}$ produces the two-fold covering of Lie groups $\pi: \SU(2) \to \SO(3)$. Now, there are two possibilities for a finite subgroup $G \subset \SU(2)$. Either we have an induced group extension:
	$$1 \to \mathbb{Z}/2\mathbb{Z} \to G \to \pi(G) \to 1,$$
	This occurs if and only if $G$ is of even order, or we have an isomorphism $G \cong \pi(G)$, if and only if $G$ is a cyclic group of odd order. Therefore, the classification of finite subgroups of $\SL(2, \mathbb{C})$ can be reduced to the classification of finite subgroups of $\SO(3)$. The problem of classifying finite rotation groups in three-dimensional Euclidean space dates back to the ancient Greeks. Indeed, a finite subgroup of $\SO(3)$ is either isomorphic to a cyclic group, a dihedral group, or the group of symmetries of a Platonic solid (see, for instance, Chapter 8 of \cite{kostrikin1982introduction} for a proof).
	
	\begin{table}[H]
		\centering
		\begin{tabular}{|c|c|c|c|}
			\hline
			Type & Name & Group & Equation \\
			\hline
			$A_n$ & Cyclic & $\mathbb{Z}/(n+1)\mathbb{Z}$ & $x^{n+1} + y^2 + z^2 = 0$ \\
			\hline
			$D_n$ & Dihedral & $D_{2(n-2)}$ & $x^{n-1} + xy^2 + z^2 = 0$ \\
			\hline
			$E_6$ & Tetrahedral & $A_4$ & $x^4 + y^3 + z^2 = 0$ \\
			\hline
			$E_7$ & Octahedral & $S_4$ & $x^3y + y^3 + z^2 = 0$ \\
			\hline
			$E_8$ & Icosahedral & $A_5$ & $x^5 + y^3 + z^2 = 0$ \\
			\hline
		\end{tabular}
		\caption{The finite subgroups of $\SO(3)$ and the corresponding Kleinian surface singularities.}
		\label{tab:Kleinian}
	\end{table}
	
	Furthermore, Klein completely solved the invariant theory for each finite subgroup $G \subset \SL(2, \mathbb{C})$. Specifically, $G$ acts on the polynomial algebra $\mathbb{C}[u, v]$ by means of linear changes of variables. As Klein observed, there exists a triple of invariant polynomials $x, y, z \in \mathbb{C}[u, v]^G$ subject to a single algebraic relation $f(x, y, z) = 0$ for some quasi-homogeneous polynomial $f \in \mathbb{C}[x, y, z]$ such that $\mathbb{C}[u, v]^G = \mathbb{C}[x, y, z]$. This implies that the algebra of invariant polynomials $\mathbb{C}[u, v]^G$ is finitely generated. Therefore, the corresponding quotient variety $X = \mathbb{A}^2_{\mathbb{C}}/G := \Spec(\mathbb{C}[u, v]^G)$ can be realized geometrically by means of the hypersurface $X \cong f^{-1}(0)$ defined by the regular map $f: \mathbb{A}^3_{\mathbb{C}} \to \mathbb{A}^1_{\mathbb{C}}$. The polynomials $f$ identified by Klein are explicitly given in Table \ref{tab:Kleinian}.
	
	\subsection{Du Val's Resolution of Kleinian Singularities}
	
\begin{figure}[H]
		
		\begin{tikzpicture}[scale=.4]
			\draw (-1,0) node[anchor=east]  {$A_n$};
			\foreach \x in {0,...,5}
			\draw[xshift=\x cm,thick] (\x cm,0) circle (.3cm);
			\draw[dotted,thick] (0.3 cm,0) -- +(1.4 cm,0);
			\foreach \y in {1.15,...,4.15}
			\draw[xshift=\y cm,thick] (\y cm,0) -- +(1.4 cm,0);
			
		\end{tikzpicture}
		
		\begin{tikzpicture}[scale=.3]
			\draw (-1,0) node[anchor=east]  {$D_n$};
			\foreach \x in {0,...,4}
			\draw[xshift=\x cm,thick] (\x cm,0) circle (.3cm);
			\draw[xshift=8 cm,thick] (30: 17 mm) circle (.3cm);
			\draw[xshift=8 cm,thick] (-30: 17 mm) circle (.3cm);
			\draw[dotted,thick] (0.3 cm,0) -- +(1.4 cm,0);
			\foreach \y in {1.15,...,3.15}
			\draw[xshift=\y cm,thick] (\y cm,0) -- +(1.4 cm,0);
			\draw[xshift=8 cm,thick] (30: 3 mm) -- (30: 14 mm);
			\draw[xshift=8 cm,thick] (-30: 3 mm) -- (-30: 14 mm);
		\end{tikzpicture}
		
		\begin{tikzpicture}[scale=.4]
			\draw (-1,1) node[anchor=east]  {$E_6$};
			\foreach \x in {0,...,4}
			\draw[thick,xshift=\x cm] (\x cm,0) circle (3 mm);
			\foreach \y in {0,...,3}
			\draw[thick,xshift=\y cm] (\y cm,0) ++(.3 cm, 0) -- +(14 mm,0);
			\draw[thick] (4 cm,2 cm) circle (3 mm);
			\draw[thick] (4 cm, 3mm) -- +(0, 1.4 cm);
		\end{tikzpicture}
		
		\begin{tikzpicture}[scale=.4]
			\draw (-1,1) node[anchor=east]  {$E_7$};
			\foreach \x in {0,...,5}
			\draw[thick,xshift=\x cm] (\x cm,0) circle (3 mm);
			\foreach \y in {0,...,4}
			\draw[thick,xshift=\y cm] (\y cm,0) ++(.3 cm, 0) -- +(14 mm,0);
			\draw[thick] (4 cm,2 cm) circle (3 mm);
			\draw[thick] (4 cm, 3mm) -- +(0, 1.4 cm);
		\end{tikzpicture}
		
		\begin{tikzpicture}[scale=.4]
			\draw (-1,1) node[anchor=east]  {$E_8$};
			\foreach \x in {0,...,6}
			\draw[thick,xshift=\x cm] (\x cm,0) circle (3 mm);
			\foreach \y in {0,...,5}
			\draw[thick,xshift=\y cm] (\y cm,0) ++(.3 cm, 0) -- +(14 mm,0);
			\draw[thick] (4 cm,2 cm) circle (3 mm);
			\draw[thick] (4 cm, 3mm) -- +(0, 1.4 cm);
		\end{tikzpicture}
		
		\caption{The simply laced Dynkin diagrams}
		\label{fig:Dynkin}
	\end{figure}
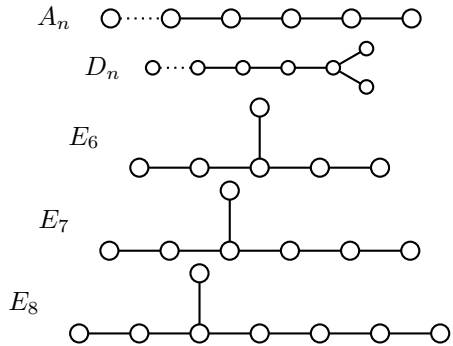

	The simply laced Dynkin diagrams (see Fig. \ref{fig:Dynkin}) are at the core of the so-called $\mathrm{ADE}$ pattern, which governs several distinct classification schemes in mathematics. For instance, simple Lie algebras, Coxeter reflection groups, quivers of finite type, and the Platonic solids are all primary examples of this mysterious universal pattern (see \cite{hazewinkel1977ubiquity} for an exposition on this fascinating subject).
	
	In the 1920s, P. Du Val \cite{du1934isolated} discovered a remarkable relationship between Kleinian surface singularities and the simply laced Dynkin diagrams, unveiling another instance of the $\mathrm{ADE}$ pattern. The idea is as follows: for any finite subgroup $G \subset \SL(2, \mathbb{C})$, the quotient variety $X := \mathbb{A}^2_{\mathbb{C}} / G$ has a single isolated singularity that can be resolved by a finite sequence of blow-ups at certain points. In fact, we can produce a minimal resolution $\pi: Y \to X := \mathbb{A}^2_{\mathbb{C}} / G$ in the sense of birational geometry. Each irreducible component of the corresponding exceptional divisor $E \subset Y$ is a rational $(-2)$-curve, which means that it is isomorphic to the projective line $\mathbb{P}^1_{\mathbb{C}}$ with self-intersection $-2$ within the smooth surface $Y$. Interestingly, the dual graph of the exceptional divisor $E$, which has a vertex for each of its irreducible components with an edge connecting them whenever they intersect, coincides with the corresponding simply laced Dynkin diagram. The type of this diagram is the same as that corresponding to the group $G$ according to the $\mathrm{ADE}$ classification in Table \ref{tab:Kleinian}.

	\begin{figure}[H]
		\begin{picture}(80,70)
			\put(-87,35){$E_1$}
			\put(-90,10){\line(1,1){50}}
			\put(-27,35){$E_2$}
			\put(-60,60){\line(1,-1){50}}
			\put(5,30){$\cdots$}
			\put(30,10){\line(1,1){50}}
			\put(66,22){$E_{n-1}$}
			\put(60,60){\line(1,-1){50}}
			\put(90,10){\line(1,1){50}}
			\put(170,0){$E_{n+1}$}
			\put(128,60){\line(1,-1){45}}
		\end{picture}
		\caption[CrepantAn]{The exceptional divisor in the crepant resolution of the surface $A_n$-singularity.}
		\label{fig:Crepant}
	\end{figure}
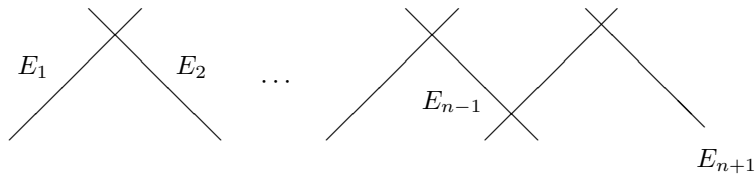
	
	Moreover, Du Val also provided a characterization of Kleinian surface singularities as those isolated surface singularities that do not affect the conditions of adjunction (see \cite{durfee1979fifteen} for a discussion of the various ways to characterize Kleinian surface singularities).
	
	\subsection{The Classical McKay Correspondence}
	
	As the 20th century progressed, several distinct aspects of Kleinian surface singularities and their role in the $\mathrm{ADE}$ classification became better understood. Notably, Brieskorn's solution \cite{brieskorn1970singular} in 1970 of a conjecture by A. Grothendieck demonstrated the realization of Kleinian surface singularities in terms of the nilpotent variety of the simply laced semi-simple Lie algebra of the corresponding type. However, a straightforward relationship between the simply laced Dynkin diagrams and binary polyhedral groups was discovered by J. McKay \cite{mckay1983graphs} in 1980 within the realm of representation theory. As McKay observed, for any finite subgroup $G$ of $\SL(2, \mathbb{C})$, we can associate a certain quiver $\widetilde{Q}_G$ in the following manner. Let $V$ denote the standard representation of $G$ given by the inclusion $G \subset \SL(2, \mathbb{C})$. The vertices of the quiver $\widetilde{Q}_G$ correspond to the irreducible representations of $G$, and there is an edge connecting a pair of irreducible representations $V_{\alpha}$ and $V_{\beta}$ whenever there is an inclusion $V_{\alpha} \subset V \otimes V_{\beta}$. The resulting quiver $\widetilde{Q}_G$ is often called the extended McKay graph of $G$. On the other hand, the sub-quiver $Q_G \subset \widetilde{Q}_G$, obtained by deleting the vertex corresponding to the trivial representation of $G$, is referred to as the McKay quiver of $G$.

	\begin{figure}[H]
		
		\begin{tikzpicture}[scale=0.7]
			\draw[xshift=-1 cm,thick] (30: 17 mm) circle (.3cm);
			\draw[xshift=-1 cm,thick] (-30: 17 mm) circle (.3cm);
			\foreach \x in {1,2,3}
			\draw[xshift=\x cm,thick] (\x cm,0) circle (.3cm);
			\draw[xshift=6 cm,thick] (30: 17 mm) circle (.3cm);
			\draw[xshift=6 cm,thick] (-30: 17 mm) circle (.3cm);
			
			\draw[xshift=-1 cm,thick] (2.5: 27.5 mm) -- (21.5: 18.5 mm);
			\draw[xshift=-1 cm,thick] (-2.5: 27.5 mm) -- (-21.5: 19 mm);
			\foreach \y in {1.15,...,2.15}
			\draw[xshift=\y cm,thick] (\y cm,0) -- +(1.4 cm,0);
			\draw[xshift=6 cm,thick] (30: 3 mm) -- (30: 14 mm);
			\draw[xshift=6 cm,thick] (-30: 3 mm) -- (-30: 14 mm);
		\end{tikzpicture}
		
		\caption{The extended Mckay graph of the dihedral group $D_3$.}
		\label{fig:MckayQuiver}
	\end{figure}
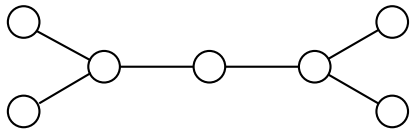
	
	In these terms, the McKay correspondence can be stated as an isomorphism of quivers:
	$$Q_G \cong E^\vee,$$
	where $E^\vee$ denotes the quiver given by the dual graph of the exceptional divisor $E$ of the minimal resolution $\pi: Y \to X := \mathbb{A}^2_{\mathbb{C}} / G$. This implies a one-to-one correspondence between the non-trivial representations of $G$ and the irreducible components of the exceptional divisor $E$. A more geometric perspective on the above isomorphism is due to González-Sprinberg and Verdier \cite{gonzalez1983construction}, who associated to each irreducible representation $V_{\alpha}$ of $G$ a vector bundle $\rho_{\alpha}$ over $Y$ such that the Chern classes $c_1(\rho_{\alpha})$ provide a basis for the Picard group $\mathsf{Pic}(Y)$. This construction can be extended to an isomorphism in $K$-theory:
	$$\Rep(G) \cong K_G(\mathbb{A}^2_{\mathbb{C}}) \cong K(Y).$$
	Moreover, a stronger version of McKay's correspondence states that there is a triangulated equivalence of derived categories (see \cite{bridgeland2001mckay}):
	$$D^b_G(\mathbb{A}^2_{\mathbb{C}}) \simeq D^b(Y).$$
	This equivalence can be interpreted as the ability to read off the geometry of the minimal resolution $Y$ directly from the representation theory of the finite group $G$, and vice versa.
	
	\subsection{The Spectrum of an Isolated Hypersurface Singularity}
	
	The notion of the spectrum of isolated hypersurface singularities was introduced in 1976 by J. Steenbrink \cite{steenbrink1976limits} based on previous work by V. Arnold \cite{arnold1981some} (we also refer the reader to \cite{van2020spectrum} for an introduction to this subject). Let us briefly recall some relevant notions. Consider a polynomial $f \in \mathbb{C}[x_0, \ldots, x_n]$ with $f(0) = 0$ such that the central fiber $f^{-1}(0) \subset \mathbb{A}^n_{\mathbb{C}}$ has an isolated singularity at the origin. A classical result by J. Milnor \cite{milnor2016singular} states that there exists a ball $X := B_{\delta} \subset \mathbb{A}^{n+1}_{\mathbb{C}}$ as well as a disc $S := \Delta_{\epsilon} \subset \mathbb{A}^1_{\mathbb{C}}$ such that the restriction $f: X^* \to S^*$ is a $C^{\infty}$-smooth fibration, where $S^* := \Delta_{\epsilon} - \{0\}$ and $X^* := X - X_0$ for $X_0 := f^{-1}(0) \cap B_{\delta}$. In particular, the fibers $X_t := f^{-1}(t) \subset X$ for any value of the parameter $t \in S^*$ are all diffeomorphic. Moreover, each $X_t$ has the same homotopy type as a bouquet of $n$-spheres, as many as the dimension of the Jacobian algebra of $f$:
	$$\mu := \dim_{\mathbb{C}} \mathbb{C}[x_0, \ldots, x_n] / (\partial_0 f, \ldots, \partial_n f).$$
	Thus, by circling the origin $0 \in S$ along a path starting and ending at any fixed value $t$, we get the corresponding monodromy operator:
	$$T: H^n(X_t) \to H^n(X_t).$$

	\begin{figure}[H]
		\begin{tikzpicture}
			\draw (2,2) circle (2cm);
			\draw (2,-2) ellipse (1.5cm and 0.5cm);
			\filldraw[black] (2,-2) circle (1pt) node[anchor=west]{0};
			\filldraw[black] (1.5,-2) circle (1pt) node[anchor=west]{t};
			\filldraw[black] (2,2) circle (1pt) node[anchor=west]{};
			\node (X0) at (3.2, 4.2) {$X_0$};
			\node (Xt) at (2.2, 4.2) {$X_t$};
			\node (X) at (2, 0) {};
			\node (S) at (2, -1.5) {};
			\node (f) at (2.3, -0.7) {$f$};
			\draw [->] (X) -- (S);
			\draw (2,2) to [out=10,in=-80] (3,4);
			\draw (2,2) to [out=10, in=80] (3,0);
			
   \draw plot [smooth, tension=2] coordinates { (2.5,-0.2) (2.3,1) (1,2) (2.3,3) (2.5,4.2)};
		\end{tikzpicture}
		\caption{The Milnor fibration.}
	\end{figure}
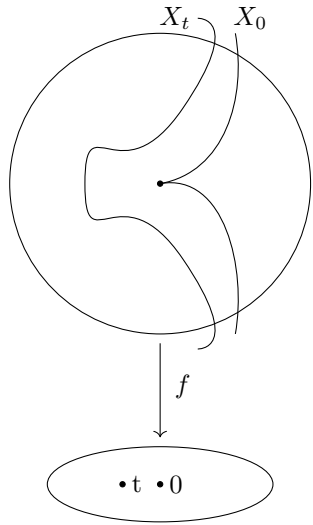
	
	The celebrated Brieskorn's monodromy theorem \cite{brieskorn1970monodromie} states that the above-defined operator $T$ is quasi-unipotent. More precisely, we have $(T^q - 1)^{n+1} = 0$ for some positive integer $q$. Hence, the eigenvalues of $T$ are all roots of unity. On the other hand, a result by Steenbrink \cite{steenbrink1974mixed} states that there is a Mixed Hodge structure on the cohomology $H^n(X_t)$. With the aid of this mixed Hodge structure, we can choose specific logarithms $\lambda$ of the eigenvalues $e^{2\pi\lambda\sqrt{-1}}$ of $T$ together with some multiplicity $m_\lambda$. These rational numbers $\lambda$ are called the spectral numbers of the singularity and can be conveniently organized into the so-called spectral polynomial:
	$$\mathrm{Sp}(f) = \sum_{\lambda \in \mathbb{Q}} m_\lambda t^{\lambda} \in \mathbb{Z}[t^\lambda, \lambda \in \mathbb{Q}].$$
	The spectrum of an isolated hypersurface singularity is a powerful invariant due to its remarkable properties. For instance, the semi-continuity of the spectrum under deformations led to a proof by A. Varchenko \cite{varchenko1983semicontinuity} of a conjecture by V. Arnold regarding an upper bound for the number of ordinary double points in a projective hypersurface.
	
	\begin{table}[H]
		\begin{tabular}{|c|c|}\hline
			Type & Spectrum \\ \hline
			$A_n$ & $\{\frac{1}{n+1},\frac{2}{n+1},...,\frac{n}{n+1}\}$ \\ \hline
			$D_n$ & $\{\frac{1}{2n-2},\frac{3}{2n-2},...,\frac{2n-3}{2n-2},\frac{n-2}{2n-2}\}$ \\ \hline
			$E_6$ & $\{\frac{1}{12},\frac{4}{12},\frac{5}{12},\frac{7}{12},\frac{8}{12},\frac{11}{12}\}$ \\ \hline
			$E_7$ & $\{\frac{1}{18},\frac{5}{18},\frac{7}{18},\frac{9}{18},\frac{11}{18},\frac{13}{18},\frac{17}{18}\}$ \\ \hline
			$E_8$ & $\{\frac{1}{30},\frac{7}{30},\frac{11}{30},\frac{13}{30},\frac{17}{30},\frac{19}{30},\frac{23}{30},\frac{29}{30}\}$ \\ \hline
		\end{tabular}
		\caption{The spectrum of the Kleinian surface singularities.}\label{tab:Spectrum}
	\end{table}
	For instance, the spectrum of the Kleinian surface singularities in Table \ref{tab:Kleinian} is depicted in Table \ref{tab:Spectrum}. Naturally, a question arises: \emph{Is there a purely representation-theoretic interpretation for the spectral numbers associated with a finite subgroup of $\SL(2, \mathbb{C})$?}
	
	\subsection{A Look Through the Mirror}
	
	Since the spectacular achievement by Candelas, de la Ossa, Green, and Parkes \cite{candelas1991pair} regarding the counting of rational curves on a quintic threefold, the phenomenon of mirror symmetry observed in string theory has evolved into a rich theory at the interface between physics and mathematics. From a mathematical perspective, a possible rigorous formulation of this phenomenon is given by the celebrated homological mirror symmetry conjecture proposed by M. Kontsevich \cite{kontsevich1995homological} at the International Congress of Mathematicians in 1994. In broad terms, certain types of geometric entities in a wide sense (ranging from varieties, stacks, singularities, or noncommutative spaces, etc.) can be associated with two types of models: the $A$-model (related to symplectic geometry) and the $B$-model (related to complex geometry), encoded by suitably defined categories of $A$-branes and $B$-branes, respectively. Homological mirror symmetry predicts that for any such geometric entity $\mathcal{X}$, there exists a mirror dual $\mathcal{X}^{\vee}$ such that the corresponding $A$ and $B$ models are interchanged by means of categorical equivalences (we refer the reader to \cite{bocklandt2021gentle} for an accessible introduction to homological mirror symmetry):
	$$A(\XX)\simeq B(\XX^{\vee})\text{ and }B(\XX)\simeq A(\XX^{\vee}).$$
	
	In the context of the present work, homological mirror symmetry reveals itself in the following manner. For any finite subgroup $G$ of $\SL(2, \mathbb{C})$, let us consider the corresponding quotient orbifold $\mathcal{X} := \mathbb{A}^2_{\mathbb{C}} \sslash G$. It turns out that the mirror dual of $\mathcal{X}$ is given by the polynomial map $f: \mathbb{A}^2_{\mathbb{C}} \to \mathbb{A}^1_{\mathbb{C}}$ associated with $G$ as depicted in Table \ref{tab:Kleinian}. In physical jargon, this map is the superpotential defining a Landau-Ginzburg model. As we will explain in Section \ref{section:Aside}, the $A$-model attached to the polynomial $f: \mathbb{A}^2_{\mathbb{C}} \to \mathbb{A}^1_{\mathbb{C}}$ can be realized using the category of twisted complexes of the so-called directed Fukaya category $\Fuk^{\to}(f)$, which can be thought of as a categorification of a distinguished basis of vanishing cycles of the Milnor fibration $f: X^* \to S^*$, that is, those cycles in $H^n(X_t)$ that collapse to zero as $t \to 0$. 
	
	On the other hand, as we will explain in Section \ref{section:Bside}, the $B$-model attached to $\mathcal{X}$ is based on the category of perfect complexes $\Perf(\mathcal{X}) \cong \Perf_G(\mathbb{A}^2_{\mathbb{C}})$, that is, bounded complexes of $G$-equivariant vector bundles over the complex affine plane $\mathbb{A}^2_{\mathbb{C}}$. Hence, as we argue in Section \ref{section:MckayMirror}, the classical McKay correspondence can be stated as a triangulated equivalence:
	$$\Tw(\Fuk^{\to}(f)) \simeq \Perf_G(\mathbb{A}^2_{\mathbb{C}}).$$
	
	Of course, the above derived equivalence provides only half of the whole homological mirror symmetry scheme taking place here. Let us mention that the $B$-model associated with the map $f: \mathbb{A}^2_{\mathbb{C}} \to \mathbb{A}^1_{\mathbb{C}}$ can be related to the so-called category of matrix factorizations $\MF(f)$, which turns out to be equivalent to the category of singularities $\Singu(X_0)$ over the central fiber $X_0 := f^{-1}(0)$ (see \cite{orlov2003triangulated}). On the other hand, the $A$-model associated with the orbifold $\mathcal{X}$ should be given by a suitably defined equivariant Fukaya category $\Fuk_G(\mathbb{A}^2_{\mathbb{C}})$. Thus, the remaining half of the homological mirror symmetry scheme would be determined by a corresponding triangulated equivalence:
	$$\MF(f) \simeq \Tw(\Fuk_G(\mathbb{A}^2_{\mathbb{C}})).$$
	
	Even though it would be quite interesting to have a complete picture of the above equivalence, as well as an interpretation of its meaning, we do not attempt to pursue these problems in the present work.

	\subsection{The Stringy Spectrum of Orbifolds}
	
	As the reader probably guessed, the directed Fukaya category $\Fuk^\to(f)$ encodes a great deal of the singularity invariants of $f$ in some way. In fact, the spectrum of $f$ can be extracted from the category of twisted complexes $\Tw(\Fuk^\to(f))$, as we will explain in Section \ref{section:Aside}. On the other hand, due to the triangulated equivalence $\Tw(\Fuk^\to(f)) \simeq \Perf_G(\mathbb{A}^2_{\mathbb{C}})$ given by the McKay correspondence, it is reasonable to ask how one can compute the spectrum of $f$ in terms of $\Perf_G(\mathbb{A}^2_{\mathbb{C}})$. Moreover, it turns out that the latter category is closely related to the representation theory of $G$, as we will explain in Section \ref{section:Bside}. Hence, one should be able to define the spectrum of $f$ in terms of the representation-theoretic data coming from $G$. Indeed, as we will explain in Section \ref{section:Spectrum}, one can define the stringy spectrum of the orbifold $\mathcal{X} = \mathbb{A}^2_{\mathbb{C}} \sslash G$ as an invariant of the category $\Perf_G(\mathbb{A}^2_{\mathbb{C}})$, and this coincides with the classical spectrum of $f$, answering the question posed at the end of Section 1.4.
	
	\subsection{Final Considerations}
	
	We would like to wrap up this introduction by pointing out that there are multiple formulations of the mirror duality for $\mathrm{ADE}$ singularities in the existing literature, with various interpretations and levels of generality. For example, but not limited to, you can refer to \cite{batyrev1996strong,kajiura2007matrix,milnor2016singular,ebeling2016homological,fan2013witten}. As far as we are aware, the specific formulation of the classical McKay correspondence as a mirror duality that we present here is novel. However, it's important to note that as this is a survey article, we haven't included many details, as those will be included in future works. Additionally, the connection we've described in the last section of this article, regarding the representation of the spectrum in terms of orbifold data, is a new perspective. It provides the initial evidence that motivated us to develop a much more comprehensive theory of the stringy spectrum of orbifolds. This work serves as an initial announcement of the results, which will be further detailed in \cite{BecerrLupercioKatzarkovSpectrum}.
	
	\section{The A-side: The category of vanishing cycles}\label{section:Aside}
	
	In this section, we will briefly describe the $A$-model associated with a Kleinian surface singularity. Fix a polynomial map $f:\aff^2_\complex\to\aff^1_\complex$ as given in Table \ref{tab:Kleinian} from now on. The last map can be thought as the superpotential defining a Landau-Ginzburg model. In this setting, one can construct the so-called directed Fukaya category $\Fuk^\to(f)$, an $A_\infty$-category that enhances the collection of vanishing cycles of $f$. In this sense, the category $\Fuk^\to(f)$ acts as a receptacle for the singularity invariants attached to the surface defined by the zero loci $V(f):=f^{-1}(0)\subset\aff^2_\complex$. Moreover, the corresponding category of complex twisted $\Tw(\Fuk^\to(f))$ with its triangulated structure provides a suitable model for the $A$-branes of the Landau-Ginzburg model defined by $f$ (we refer the reader to Chapter 3 of \cite{bocklandt2021gentle} for a readable discussion on this notion).    
	\subsection{The Milnor fibration} Fortunately, the polynomials of Table \ref{tab:Kleinian} belong to a class of polynomials for which the corresponding singularity invariants can be described in rather explicit manner; namely, the class of weighted homogeneous polynomials. Consider a polynomial $$f=\sum_{n=(n_x,n_y,n_z)\in\naturals^3}c_n x^{n_x}y^{n_y}z^{n_z}$$ in $\complex[x,y,z]$. The support of $f$ is defined as the set $\{n\in\naturals^3|c_n\not=0\}$ of $\naturals^3$. We say that $f$ in $\complex[x,y,z]$ is weighted homogeneous with weights $(w_x,w_y,w_z)\in\rationals_{>0}$ if its support is contained in the hyperplane of with equation $$w_xn_x+w_yn_y+w_zn_z=1.$$ Thus, if $d$ the smallest common denominator of $w_x,w_y,w_z$ and we define the integer numbers $a=dw_x$, $b=dx_y$ and $c=dw_z$, we have $$f(\lambda^ax,\lambda^by,\lambda^cz)=\lambda^df(x,y,z)$$ for $\lambda\in\complex^*$. In such a case, the restriction $f:\aff^3_{\complex}-f^{-1}(0)\to \aff^1_{\complex}-\{0\}$ turns out to be a $C^{\infty}$-smooth fibration, called the global affine Milnor fibration of $f$. By moving around the unit circle in $\aff^1_{\complex}-\{0\}$, we obtain the corresponding monodormy homeomorphism $h:F\to F$ where $F:=f^{-1}(1)\subset\aff^3_{\complex}$ is the global affine Milnor fiber and $$h(x,y,z):=e^{\frac{2\pi\sqrt{-1}}{d}}\cdot(x,y,z).$$ Moreover, if $f$ has just an isolated singularity at the origin $0\in\aff^3_{\complex}$, then the associate Milnor number is given by $$\mu:=\dim\Jac(f)=\frac{d-a}{a}\cdot\frac{d-b}{b}\cdot\frac{d-c}{c}$$ where $$\Jac(f):=\complex[x,y,z]/(\partial_xf,\partial_yf,\partial_zf)$$ is the Jacobian algebra of $f$. On the other hand, it is well-known that the global affine Milnor fiber $F$ has the homotopy type of a bouquet of spheres $\bigvee_{i=1}^\mu S^2$, so that the cohomology of $F$ is given by: $$H^n(F,\complex)=\begin{cases} \complex & \text{if n=0} \\
		\complex^{\mu} & \text{if n=2} \\
		0 & \text{otherwise}\end{cases}$$ In fact, we have the following (cf. Theorem 9.2.8 of \cite{steenbrink2022mixed}):
	
	\begin{theorem} There are identifications $H^2(F,\complex)\cong H^2_\DR(F)\otimes\complex\cong\Omega_f,$ where $$\Omega_f:=\Omega^3_{\aff^3_\complex}/df\wedge\Omega^2_{\aff^3_\complex}\cong\Jac(f)\cdot dx\wedge dy\wedge dz$$ is the so-called Milnor module of $f$, induced by the residue map $$\omega\in\Omega^3_{\aff^3_\complex}\mapsto\Res_F(\omega)\in\Omega^2_F$$ where $\eta:=\Res_F(\omega)$ satisfies $\omega=\eta\wedge\frac{df}{f}+\omega'$ for some regular $3$-form $\omega'$ in $\aff^3_\complex$ \end{theorem}

	\subsection{Vanishing cycles} In fact, we can do better. Let $A\subset\integers^3_{\geq 0}$ be a finite sub-set such that the monomials $x^{\alpha_x}y^{\alpha_y}z^{\alpha_z}$ for $\alpha=(\alpha_x,\alpha_y,\alpha_z)\in A$ for a basis of the Jacobian algebra $\Jac(f)$. For $\alpha$ in $A$, define $$n(\alpha):=\omega_x(\alpha_x+1)+\omega_y(\alpha_y+1)+\omega_z(\alpha_z+1)$$ and $$\omega_\alpha:=\frac{x^{\alpha_x}y^{\alpha_y}z^{\alpha_z}}{(f-1)^{[n(\alpha)]}}dx\wedge dy\wedge dz.$$ Thus, we have the following (which is a particular case of Theorem 9.2.8 of \cite{steenbrink2022mixed}):
	
	\begin{theorem} The collection $\{[\eta_\alpha]\in H^2(F,\complex)\}_{\alpha\in A}$ where $\eta_\alpha:=\Res_F(\omega_\alpha)$ form a basis of $H^2(F,\complex)$.\end{theorem}
	
	Moreover, one can also prove a geometric version of the above cohomological result:
	
	\begin{theorem}\label{thm:VanishingCycles} For each $\alpha$ in $A$, let $V_\alpha\subset F$ be the Poincar\'e dual of $\eta_\alpha$. Then, the collection $\{V_\alpha\}_{\alpha\in A}$ is distiguished basis of vanishing cycles for $f$.\end{theorem}
	
	\subsection{The directed Fukaya category} We define $\Fuk^\to(f)$ as the $A_\infty$-category whose objects are given by a distinguished basis of vanishing cycles $\{V_\alpha\}_{\alpha\in A}$ as in the last paragraph and whose morphism spaces are given by Lagrangian intersection Floer complexes. Roughly speaking, the corresponding triangulated category of twisted complexes $\Tw(\Fuk^\to(f))$ has for objects the pairs $(C,d)$ where $C$ belongs to the additive completion of $\Fuk^\to(F)$ and $d$ is a cohomological differential operator on $C$ satisfying the Maurer-Cartan equation (see Chapter 3 of \cite{bocklandt2021gentle} for details) 
	
	\section{The B-side: The category of perfect complexes}\label{section:Bside}
	
	The goal of this section is to provide a brief description of the $B$-model associated to a Kleinian quotient singularity. For any finite sub-group $G$ of $\SL(2,\complex)$, the corresponding quotient variety $X:=\aff^2_\complex/G$ is affine $X\cong\Spec(\complex[x,y]^G)$ with a single isolated singularity. On the other hand, the natural $G$-action on the canonical bundle $\omega_{\aff^2_\complex}$ is trivial and, therefore, it descends into the canonical sheaf $\omega_X$, which is trivial $\omega_X\cong\OO_X$. In other words, the variety $X$ is a singular and non-compact Calabi-Yau surface. In fact, we can do better by considering instead the quotient stack $\XX:=\aff^2_\complex\sslash G$, which is a Calabi-Yau orbifold. As is explain in Section \ref{section:MckayMirror}, the orbifold $\XX$ turns out to be the mirror dual of the Landau-Ginzburg model given by the corresponding Kleinian polynomial $f_G:\aff^2_\complex\to\aff^1_\complex$ as in Table \ref{tab:Kleinian}. In order to see this, we will define the $B$-model attached to $\XX$ as the triangulated category of perfect complexes $B(\XX):=\Perf(\XX)$, as we presently describe in more detail.     
	
	\subsection{The category of perfect complexes} Let us start by recalling that there is a linear equivalence $\Coh(\XX)\cong\Coh_G(\aff^2_\complex)$ between the category of coherent sheaves over $\XX$ and the category of $G$-equivariant coherent sheaves over $\aff^2_\complex$. Thus, we also have a triangulated equivalence between derived categories $D^b(\XX)\cong D^b_G(\aff^2_\complex)$. Let us denote by $\Perf(\XX)\subset D^b(\XX)$ the subcategory of perfect complexes, that is, those bounded complexes of coherent sheaves over $\XX$ that are locally quasi-isomorphic to a bounded complex of vector bundles in the \'etale topology of $\XX$. However, since the orbifold $\XX=\aff^2_\complex\sslash G$ is a global quotient, it has the resolution property (which means that every coherent sheaf over $\XX$ is the quotient of a vector bundle), and this implies in turn that every perfect complex over $\XX$ is in fact quasi-isomorphic to a bounded complex of global vector bundles over $\XX$ (see \cite{totaro2004resolution}). Let us denote by $\Vect(\XX)\subset\Coh(\XX)$ the exact sub-category of vector bundles over $\XX$ (so that we have a corresponding equivalence $\Vect(\XX)\cong\Vect_G(\aff^2_\complex)$). Thus, we have the following:
	
	\begin{theorem} There is a triangulated equivalence $\Perf(\XX)\simeq D^b(\Vect(\XX))$.\end{theorem}
	
	\subsection{The skew group algebra} There is rather nice algebraic description of the category of perfect complexes $\Perf(\XX)$ which runs as follows. Note that the abelian category $\Coh(\XX)$ has a projective generator. Indeed, if $\complex[G]$ denotes the regular representation of $G$, the trivial vector bundle $\complex[G]\times\aff^2_\complex\to\aff^2_\complex$ can be endowed with a $G$-equivariant structure by means of the diagonal $G$-action on the total space $\complex[G]\times\aff^2_\complex$. We will denote by $\EE_{\complex[G]}$ the corresponding stacky vector bundle over $\XX=\aff^2_\complex\sslash G$. Hence, it is easy to see that any coherent sheaf $\FF$ over $\XX$ admits an epimorphism $\EE_{\complex[G]}^n\to\FF\to 0$ for some positive integer $n$. Thus, $\EE_{\complex[G]}$ provides a projective generator of $\Coh(\XX)$, as have been claimed. Note that, for any coherent sheaf $\FF$ over $\XX$, the linear space $\Hom(\EE_{\complex[G]},\FF)$ acquires the structure of a module over the endomorphism algebra $\End(\EE_{\complex[G]})$ by composition of morphisms. The last immediately implies that:
	
	\begin{theorem} The assignment $\FF\mapsto\Hom(\EE_{\complex[G]},\FF)$ induces a linear equivalence of categories $\Coh(\XX)\simeq\Mod(\End(\EE_{\complex[G]}))$\end{theorem}
	
	In fact, it is well-known that there is a much more concrete description of the endomorphism algebra $\End(\EE_{\complex[G]})$. Note that the group $G$ acts on the polynomial algebra $\complex[x,y]$. In fact, a matrix $g=\begin{pmatrix} a & b \\ c & d \end{pmatrix}$ in $G\subset\SL(2,\complex)$ acts on a given polynomial $f$ in $\complex[x,y]$ by changing variables $f^g(x,y)=f(ax+by,cx+dy)$. We can define the skew group algebra $\complex[x,y]\rtimes G$ as the $\complex$-algebra generated by the pairs $(f,g)$ where $f\in\complex[x,y]$ and $g\in G$, subject to the relations $$(f_1,g_1)\cdot (f_2,g_2)=(f_1f_2^{g_1},g_1g_2).$$ It turns out that the quotient orbifold $\XX=\aff^2_\complex\sslash G$ is affine and the skew group algebra $\complex[x,y]\rtimes G$ plays the role of its algebra of regular functions (which is not commutative in general). More precisely, we have the following: 
	
	\begin{theorem} There is an isomorphism of algebras $\End(\EE_{\complex[G]})\cong\complex[x,y]\rtimes G$. In particular, there is an equivalence $\Vect(\XX)\simeq\Proj(\complex[x,y]\rtimes G)$ between the category of vector bundles over $\XX$ and the category of finitely generated projective modules over $\complex[x,y]\rtimes G$, and this induces a triangulated equivalence: $$\Perf(\XX)\simeq D^b(\Proj(\complex[x,y]\rtimes G))$$\end{theorem}
	
	\subsection{Orbifold cohomology} The triangulated category $\Perf(\XX)$ can be thought as a sort of universal invariant of the orbifold $\XX$. In principle, any additive invariant of $\XX$, like the orbifold Euler number $\chi^{\orb}(\XX)$ or the orbifold $K$-theory $K^{\orb}(\XX)$, can be extracted out of $\Perf(\XX)$ by means of a definite procedure. In this work, we are mainly interested in the so-called orbifold cohomology (sometimes also referred as the Chen-Ruan cohomology) introduced in \cite{chen2004new}. Let us recall that the orbifold cohomology of $\XX$ with complex coefficients is defined by $$H_\orb^\bullet(\XX,\complex):=H^\bullet(\Lambda\XX,\complex)$$ where $\Lambda\XX$ is the inertia orbifold of $\XX$ and the right hand side in the above equality is the usual cohomology of $\Lambda\XX$. On the other hand, due to the localization principlie for orbifold theories, we can write down the orbifold cohomology of $\XX$ as $$H_\orb^\bullet(\XX,\complex)\cong\bigoplus_{\alpha\in\pi_0(\Lambda\XX)}H^{\bullet-\age(\alpha)}(\XX_\alpha,\complex)$$ where $\XX_\alpha\subset\XX$ is the connected component corresponding to $\alpha$ and $\age(\alpha)$ denotes its age (see \cite{de2004localization}). In our particular case where $\XX=\aff^2_\complex\sslash G$ is a global quotient orbifold, each $\XX_\alpha$ is contractible and we get an identification $$H_\orb^\bullet(\XX)\cong\Rep(G)\otimes_\integers\complex$$ where $\Rep(G)$ denotes the representation ring of $G$. However, note that we have an additional grading by age in $H_\orb^\bullet(\XX,\complex)$, which will become relevant later on. Finally, let us just mention that the orbifold cohomology of $\XX$ can be recovered out of $\Perf(\XX)$ as the periodic cyclic cohomology of this last (see \cite{baranovsky2003orbifold}):  
	
	\begin{theorem}\label{thm:OrbifoldCohomolgy} There is an identification $H^\bullet_\orb(\XX,\complex)\cong \mbox{HP}^\bullet(\Perf(X))$ between the orbifold cohomology of $\XX$ and the periodic cyclic cohomology of $\Perf(\XX)$.\end{theorem}
	
	\section{The classical Mckay correspondence as a mirror duality}\label{section:MckayMirror}
	
	This section will formulate the classical Mckay correspondence as a mirror duality relating the $A$-side of a Kleinian surface singularity and the $B$-side of the quotient orbifold of the associated binary polyhedral group. To this end, we will show that each side can be realized by means of the derived category of the Dynkin quiver of the corresponding type.
	
	\subsection{The simply laced Dynkin diagrams} By a quiver we just mean a finite directed graph. A representation $V$ of a quiver $Q$ consists of finite-dimensional vector space $V_q$ for each vertex $q$ of $Q$ and a linear map $T_l:V_q\to V_p$ for each directed edge $l:q\to p$. It turns out that there is an equivalence $$\Rep(Q)\cong\Mod(P_Q)$$ between the category of representations of $Q$ and the category of finitely generated modules over the path algebra $P_Q$. We will say that $Q$ is of finite type iff there is a finite collection of indecomposable representations of $Q$. As another instance of the $\mbox{ADE}$ pattern, a celebrated result of P. Gabriel \cite{gabriel10unzerlegbare} states that (we also refer the reader to Chapter 6 of \cite{etingof2011introduction} for a proof of this fact):
	
	\begin{theorem} The quivers of finite type are given exactly by the simply laced Dynkin diagrams of $\mbox{ADE}$ type as those in Figure \ref{fig:Dynkin}.\end{theorem}
	
	\subsection{The A-side} For any polynomial $f$ as those in Table \ref{tab:Kleinian}, the distinguished basis of vanishing cycles $\{V_\alpha\}_{\alpha\in A}$ given in Theorem \ref{thm:VanishingCycles} provides a strong exceptional collection in the triangulated category $\Tw(\Fuk^\to(f))$. Moreover, the corresponding semi-orthogonal decomposition of the last category can be encoded by a Dynkin quiver $Q_f$, called the quiver of vanishing cycles of $f$, having the same $\mbox{ADE}$ type as the polynomial $f$. Thus, we have the following:
	
	\begin{theorem}\label{thm:Aside} There is a triangulated equivalence $\Tw(\Fuk^{\to}(f))\simeq D^b(Q_f)$.\end{theorem}
	
	\subsection{The B-side} On the other hand, for any finite sub-group $G$ of $\SL(2,\complex)$ we have the corresponding quotient orbifold $\XX:=\aff^2_\complex\sslash G$. Hence, it turns out that the connected components of the inertia orbifold $\Lambda\XX$, which are sometimes called the twisted sector of $\XX$, give rise to a strong exceptional collection in the triangulated category $\Perf(\XX)$. Similarly as in paragraph (4.2), this induces a semi-orthogonal decomposition of $\Perf(\XX)$ which can be encoded by the Mckay quiver $Q_G$. In this manner, we have the following:  
	
	\begin{theorem}\label{thm:Bside} There is a triangulated equivalence $\Perf_G(\aff^2_{\complex})\simeq D^b(Q_G)$.\end{theorem}
	
	\subsection{The classical Mckay correspondence} Now, suppose that $G$ is a finite sub-group of $\SL(2,\complex)$ and $f:\aff^3_\complex\to\aff^1_\complex$ the corresponding polynomial map as in the classification of Table \ref{tab:Kleinian}. Thus, we have the Mckay quiver $Q_G$ and the quiver of vanishing cycles $Q_f$, respectively, as in paragraphs (4.2) and (4.3). Hence, the classical Mckay correspondence can be stated as follows:
	
	\begin{theorem}[{\bf Classical Mckay correspondence}]\label{thm:MckayClassical} There is an isomorphism of quivers $Q_f\cong Q_G$.\end{theorem}
	
	Finally, ass a collective consequence of Theorems \ref{thm:Aside},\ref{thm:Bside} and \ref{thm:MckayClassical}, we get the following homological form of the classical Mckay correspondence:
	
	\begin{theorem}[{\bf Homological Mckay correspondence}]\label{thm:MckayMirror}There is a triangulated equivalence:$$\Tw(\Fuk^{\to}(f))\simeq\Perf_G(\aff^2_{\complex})$$\end{theorem}
	
	\section{The spectrum of Kleinian surface singularities}\label{section:Spectrum}
	
	To conclude this article, we will explain how one can extract the spectrum of a Kleinian surface singularity out of the representation theoretic data coming from the corresponding binary polyhedral group, answering the question stated at the end of paragraph (1.4).
	
	\subsection{The Picard-Lefschetz formula} Let $\{V_\alpha\}_{\alpha\in A}$ be a distinguished basis of vanishing cycles as in paragraph (2.2). We denote by $\delta_\alpha$ the reflection along the normal hyperplane to the class $[\eta_\alpha]\in H^2(F,\complex)$ with respect to the intersection form of $H^2(F,\complex)$. Thus, the cohomological monodromy operator $$T:=h^*:H^2(F,\complex)\to H^2(F,\complex)$$ is given by the Picard-Lefschetz formula:$$T=\prod_{\alpha\in A}\delta_\alpha$$ Note that, in this case, the monodromy operator $T$ is of finite order. In fact, we have $T^d=1$.
	
	\subsection{The Milnor module} Note for a weighted homogeneous polynomial $f$ with weights $(w_x,w_y,w_z)\in\rationals_{{>0}}^3$, the Milnor module $$\Omega_f:=\Omega^3_{\aff^3_\complex}/df\wedge\Omega^2_{\aff^3_\complex}$$ introduced in Paragraph 2.1 has a $\rationals$-grading defined in the following manner. We define the degree of a monomial $x^{n_x}y^{n_y}z^{n_z} dx\wedge dy\wedge dz$ as:$$\deg(x^{n_x}y^{n_y}z^{n_z} dx\wedge dy\wedge dz):=a(n_x+1)+b(n_y+1)+c(n_z+1)$$ Thus, for any rational number $\lambda$ we define $\Omega_{f,\lambda}$ as the sub-module of $\Omega_f$ whose elements are given by those polynomials whose constituent monomials have degree equal to $\lambda$. In this manner, we get a decomposition:$$\Omega_f:=\bigoplus_{\lambda\in\rationals}\Omega_{f,\lambda}$$ It turns out that one can write down the spectral polynomial of $f$ explicitly as the graded dimension of $\Omega_f$, namely:
	
	\begin{theorem}\label{thm:SpectralMilnor} The identity $\Sp(f)=\sum_{\lambda\in\rationals}\dim(\Omega_{f,\lambda})t^\lambda$ holds true.\end{theorem}
	
	On the other hand, we can recover the Milnor module $\Omega_f$ out of the triangulated category $\Tw(\Fuk^\to(f))$ as follows:
	
	\begin{theorem}\label{thm:MilnorModule} There is a canonical identification $\Omega_f\cong\mbox{HP}(\Tw(\Fuk^\to(f)))$ between the Milnor module of $f$ and the periodic cyclic cohomology of $\Tw(\Fuk^\to(f))$.\end{theorem}
	
	Additionally, one can recover the Hodge filtration $F^\bullet$ of $\Omega_f$ by taking the graded piece $F^p\Omega_f$ as the subspaces generated by the elements $\eta_\alpha$ in the distinguished basis having $n(\alpha)\leq n-p+1$.
	
	\subsection{The stringy spectrum} Of course, due to the triangulated equivalence $$\Perf_G(\aff^2_\complex)\simeq\Tw(\Fuk^\to(f))$$ 
	given by the homological Mckay correspondence as in Theorem \ref{thm:MckayMirror}, we have the following consequence of Theorem \ref{thm:OrbifoldCohomolgy}:
	
	\begin{theorem}\label{thm:MilnorOrbifold} There is an identification $\Omega_f\cong H_\orb^\bullet(\XX,\complex)$ where $\XX:=\aff^2\sslash G$.\end{theorem} 
	
	Moreover, it turns out that, under the above identification $$\Omega_f\cong H_\orb^\bullet(\XX,\complex),$$ the Picard-Lefschetz formula regarding the monodromy operator on $\Omega_f$ can be translated in terms of the Frobenius algebra structure of $H_\orb^\bullet(\XX,\complex)$, which is given using the so-called Chen-Ruan product (see \cite{chen2004new}). One can also recover the $\rationals$-grading of $\Omega_f$ from a suitably defined $\rationals$-grading on $H_\orb^\bullet(\XX,\complex)$ which is related to the grading by age. Therefore, due to Theorem \ref{thm:SpectralMilnor}, it should be possible to recover the spectral polynomial of $f$ from the representation ring $$\Rep(G)\cong H_\orb^\bullet(\XX,\complex).$$ However, we must conclude here and refer the reader to \cite{BecerrLupercioKatzarkovSpectrum} for a complete account of this story.
	
	\bibliographystyle{plain}
	\bibliography{references}
	
\end{document}